\documentclass{article}
\usepackage{amsmath}
\usepackage{amsthm}
\usepackage{amsfonts}
\usepackage{amssymb}
\usepackage{latexsym} 
\usepackage[matrix,tips,graph,curve]{xy}

\makeatletter

\def\@maketitle{%
  \newpage
  \null
  \let \footnote \thanks
    {\normalfont\sffamily\bfseries\Large\noindent\@title \par}%
    \vskip 1em%
    {\normalfont\sffamily\large
        \noindent
        \@author
        \par}
  \par
  \vskip 4em}
\def\@seccntformat#1{\csname the#1\endcsname{.\ }}
\renewcommand\section{\@startsection {section}{1}{\z@}%
                                   {-3.0ex \@plus -1ex \@minus -.2ex}%
                                   {1.5ex \@plus.2ex}%
                                   {\normalfont\large\bfseries}}
\renewcommand\subsection{\@startsection{subsection}{2}{\z@}%
                                     {-2.75ex\@plus -1ex \@minus -.2ex}%
                                     {1.5ex \@plus .2ex}%
                                   {\normalfont\large}}
\makeatother

\makeatletter 
\@addtoreset{equation}{section}
\makeatother

\theoremstyle{plain}
\newtheorem{theorem}[equation]{Theorem}
\newtheorem{corollary}[equation]{Corollary}
\newtheorem{lemma}[equation]{Lemma}
\newtheorem{proposition}[equation]{Proposition}

\theoremstyle{definition}
\newtheorem{definition}[equation]{Definition}

\newtheorem{remark}[equation]{Remark}
\newtheorem{remarks}[equation]{Remarks}

\newtheorem{example}[equation]{Example}

\newtheorem{notation}[equation]{Notation}

\newtheorem{assumption}[equation]{Assumption}
\newtheorem*{claim}{Claim}

\newcommand{\ff}{\mathfrak{f}}
\newcommand{\fg}{\mathfrak{g}}

\newcommand{\sC}{\mathcal{C}}
\newcommand{\sD}{\mathcal{D}\,}
\newcommand{\sE}{\mathcal{E}\,}
\newcommand{\sF}{\mathcal{F}\,}

\newcommand{\sI}{\mathcal{I}\,}

\newcommand{\sL}{\mathcal{L}\,}
\newcommand{\sM}{\mathcal{M}\,}
\newcommand{\sN}{\mathcal{N}\,}
\newcommand{\sO}{\mathcal{O}}
\newcommand{\sP}{\mathcal{P}\,}

\newcommand{\sR}{\mathcal{R}}
\newcommand{\sS}{\mathcal{S}}
\newcommand{\sT}{\mathcal{T}\,}

\newcommand{\IC}{\mathbb{C}}

\newcommand{\IQ}{\mathbb{Q}}
\newcommand{\IR}{\mathbb{R}}

\newcommand{\IZ}{\mathbb{Z}}

\newcommand{\End}{\mathrm{End}}

\renewcommand{\index}{\mathrm{index \,}}
\newcommand{\Hom}{\mathrm{Hom}}

\newcommand{\Spin}{\rm Spin}

\newcommand{\Ker}{\mathrm{Ker \,}}
\newcommand{\Coker}{\mathrm{Coker \,}}

\newcommand{\Hess}{\mathrm{Hess \,}}
\newcommand{\grad}{\mathrm{grad \,}}

\newcommand{\ch}{\rm ch} % Chern Character
\newcommand{\rk}{\rm rk} 

\renewcommand\dim{{\rm dim\,}}
\renewcommand\det{{\rm det\,}}

\newcommand\Hilb{\mathrm{Hilb}}
\newcommand\vol{\mathrm{vol}}

\renewcommand\i{\sqrt{-1}}

\newcommand\iso{{\cong}} 
\newcommand\tensor{{\otimes}}

\newcommand{\<}{\langle}
\renewcommand{\>}{\rangle}
\newcommand{\delbar}{\bar{\partial}}

\def\d/{/\mspace{-6.0mu}/}
\newcommand{\git}[3]{#1\d/_{\mspace{-4.0mu}#2}#3}
\newcommand\zetahilb{\zeta_{{\text{Hilb}}}}
\def\Fy{\sF \mspace{-3.0mu} \cdot \mspace{-3.0mu} y}

\def\wt{\widetilde}

\begin{document}
\title{Geometrical McKay Correspondence for Isolated Singularities}
\author{Anda Degeratu%
\footnote{Duke University; e-mail:anda@math.duke.edu}}
\date{\today}
\maketitle
%
%\begin{abstract}
% A Calabi-Yau orbifold is locally modeled on C^n/G where G is a
% finite subgroup of SL(n, C). In dimension n=3 a crepant
% resolution is given by Nakamura's G-Hilbert scheme. This crepant
% resolution has a description as a GIT/symplectic quotient. We use
% tools from global analysis to give a geometrical generalization
% of the McKay Correspondence to this case.
%\end{abstract}
%
\section*{Introduction}\label{sec:introduction}
\subsection{Crepant Resolutions of Singularities of Calabi-Yau
  orbifolds}
A Calabi-Yau manifold is a complex K\"ahler manifold with trivial
canonical bundle. In the attempt to construct such manifolds 
%it has been proved that 
%insufficient to limit our attention to smooth Calabi-Yau
%manifolds; 
it is useful to take into consideration singular Calabi-Yaus.
One of the simplest  singularities which can arise is an orbifold
singularity. An orbifold is the quotient of a smooth Calabi-Yau manifold
by a discrete group action which generically has fixed points.
Locally such an orbifold is modeled on $\IC^n/G$, where $G$ is a
finite subgroup of $SL(n, \IC)$. 

From a geometrical perspective we can try to resolve the
orbifold singularity. A resolution $(X, \pi)$ of $\IC^n/G$ is
a nonsingular complex manifold
$X$ of dimension $n$ with a proper biholomorphic map 
$\pi: X \to \IC^n/G$ that induces a biholomorphism between dense open
sets. We call $X$ a {\em crepant resolution}\footnote{Etymology: For a resolution 
of singularities we can define a notion of {\em discrepancy}~\cite{YPG}. 
A crepant resolution is a resolution without discrepancy.}   
if the canonical bundles are isomorphic, 
$K_X \iso \, \pi^*(K_{\IC^n/G})$. 
Since Calabi-Yau manifolds have trivial canonical bundle, 
to obtain a Calabi-Yau structure on $X$ one must
choose a crepant resolutions of singularities. 

It turns out that the amount of information we know about a crepant
resolution of singularities of $\IC^n/G$ depends dramatically on the
dimension $n$ of the orbifold:
\begin{enumerate}
 \item[$n=2$:] A crepant resolution always exists and is
    unique. Its topology is entirely described in terms of the finite
    group $G$ (via the McKay Correspondence).
  \item[$n=3$:] A crepant resolution always exists but it is not
    unique; they are related by flops. However all the crepant
    resolutions have the same Euler and Betti numbers: the {\em
    stringy} Betti and Hodge numbers of the orbifold~\cite{dhvw}.
  \item[$n \geq 4$:] In this case very little is known; crepant resolutions
    exist in rather special cases. Many singularities are 
    terminal, which implies that they admit no crepant resolution.
\end{enumerate}
We would like to completely understand the topology of crepant
resolutions in the case $n=3$. In this paper we are concerned with the
study of the ring structure in cohomology. This is
related to the generalization of the McKay Correspondence. In what
follows we give a description of the problem by moving
back and forth between the case $n=2$ and $n=3$.

\paragraph{The case $n=2$.}
The quotient singularities $\IC^2/G$, for $G$ a finite
subgroup of $SL(2, \IC)$, were first classified by Klein in 1884 and are
called {\em Kleinian singularities} (they are also known as  {\em Du
Val singularities} or {\em rational double points}).
There are five families of finite subgroups of $SL(2,\IC)$: 
the cyclic subgroups $\sC_k$, the binary dihedral groups
$\sD_{k}$ of order $4k$, 
the binary tetrahedral group $\sT$ of order $24$, 
the binary octahedral group $\sO$ of order $48$, 
and the binary icosahedral group $\sI$ of order $120$. 
A crepant resolution exists for each family and is unique. Moreover
the finite group completely describes the topology of the resolution. 
This is encoded in the McKay Correspondence~\cite{mckay}, which
establishes a bijection between the set of irreducible representations
of $G$ and the set of vertices of an extended Dynkin diagram of type
$ADE$ (the Dynkin diagrams corresponding to the simple Lie algebras of
the following five types: $A_{k-1}$, $D_{k+2}$, $E_6$, $E_7$ and
$E_8$).

\medskip
\noindent
Concretely, let $\{R_0, R_1, \ldots, R_r\}$ be the set of irreducible
representations of $G$, where $R_0$ denotes the one-dimensional
trivial representation.
To $G$ and its irreducible representations we associate an
$(r+1)\times(r+1)$ adjacency matrix 
$A = [a_{ij}]$ with $i, j = 0, \ldots, r$.
The entries $a_{ij}$ are positive integers; they are defined by the
tensor product decompositions
\begin{equation*}
  R_i \, \tensor \, Q = \sum_{i=0}^r \, a_{ij} \, R_{j} ,
\end{equation*}
where $Q$ denotes the natural two-dimensional representation of $G$
induced from the embedding $G \subset SL(2, \IC)$.
McKay's insight was to realize that the matrix $A$ is related to
the Cartan matrix $C$ of a Dynkin diagram of type $ADE$, via
\begin{equation}\label{eq:mckay}
  A = 2 I  - \widetilde{C}.
\end{equation}
(Here $\widetilde{C}$ is the Cartan matrix of the extended Dynkin
 diagram; the matrix $C$ is the $r\times r$-minor obtained by removing
 the first row and the first column from $\widetilde{C}$.) 

\medskip
\noindent
Using McKay's correspondence it is easy now to describe the crepant
resolution $\pi: X \to \IC^2/G$.
The exceptional divisor $\pi^{-1}(0)$ is the dual of the Dynkin
diagram: the vertices of the Dynkin diagram correspond naturally
to rational curves $C_i$ with self-intersection $-2$.
Two curves intersect transversally at one point if and only if the
corresponding vertices are joined by an edge in the Dynkin diagram,
otherwise they do not intersect. 
The curves above form a basis for $H_2(X, \IZ)$. The
intersection form with respect to this basis is the negative of the
Cartan matrix.

\medskip
\noindent
The first geometrical interpretation of the McKay Correspondence was given
by Gonzalez-Sprinberg and Verdier~\cite{GV}. To each of the
irreducible representations $R_i$ they associated a locally free
coherent sheaf $\sR_i$. The set of all these coherent sheaves form a
basis for $K(X)$, the $K$-theory of X. Moreover, the first Chern classes
$c_1(\sR_i)$ form a basis in $H^2(X, \IQ)$ and the product of 
two such classes in $H^{*} (X, \IQ)$ is given by the formula
\begin{equation}\label{eq:dim-2-multiplication}
  \left[
       \int_X  c_1 (\sR_i) c_1 (\sR_j) 
  \right]_{i,j=1, \ldots, r} = - \, C^{-1},
\end{equation}
where $C^{-1}$ is the inverse of the Cartan matrix.
The proof given by Gonzalez-Sprinberg and Verdier uses a case by case
analysis and techniques from algebraic geometry.
Kronheimer and Nakajima gave a proof of the formula
using techniques from gauge theory \cite{KN}.

\medskip
\noindent
To summarize, in the case of surface singularities,
$\IC^2/G$, the representation theory of the finite group $G$
completely determines the topology the crepant resolution. The
Dynkin diagram and the Cartan matrix (and hence the simple Lie algebra
$\fg$ associated to it) encode everything we want to know about 
the topology of the crepant resolution.

\paragraph{The case $n=3$.}
The finite subgroups of $SL(3, \IC)$ were classified by Blichfeldt in
$1917$~\cite{blichfeldt:17}: there are ten families of such finite
subgroups.  In the early $1990$'s a case by case analysis was used to
construct a crepant resolution of $\IC^3/G$ with the given stringy
Euler and Betti numbers (see~\cite{roan:96} and the references therein).
As a consequence of these constructions, we know that
all the crepant resolutions of $\IC^3/G$ have 
the Euler and Betti numbers given by the stringy Euler
and Betti numbers of the orbifold (since these numbers are unchanged under flops).
In $1995$ Nakamura made the conjecture that $\Hilb^G(\IC^3)$ is a
crepant resolution of $\IC^3/G$.
In general, for $G$ a finite subgroup of $SL(n, \IC)$, the algebraic
variety $\Hilb^G (\IC^n)$ parametrizes the $0$-dimensional
$G$-invariant subschemes of $\IC^n$ whose space of global sections is
isomorphic to the regular representation of $G$.  Nakamura made the
conjecture based on his computations for the case
$n=2$~\cite{ito-nakamura:96}; then he proved it in dimension $n=3$ for
the case of abelian groups~\cite{nakamura:99}.  In $1999$ Bridgeland,
King and Reid gave a general proof of the conjecture in the case
$n=3$, relying heavily on derived category techniques~\cite{bkr}.

\subsection{The Scope and Main Result of the paper}
In the case of surface singularities, an important feature of the
McKay Correspondence is that it gives the ring structure in
cohomology in terms of the finite group. For the case $n \geq 3$, 
nothing is known about the multiplicative structures  in cohomology or $K$-theory.
% The goal  of this project is to give an answer to the question: Can tensor
% product of $G$-modules and tensor product in $K$-theory be
% related?

We start with $X$ a crepant resolution of $\IC^n/G$. We
assume that this resolution is given by $\Hilb^{G} (\IC^n)$.
On this resolution, Ito and Nakajima, \cite{ito-nakajima:00},
gave a recipe of extending the Gonzalez-Sprinberg-Verdier sheaves, 
associating a locally free sheaf (therefore an algebraic vector
bundle) $\sR_i$ to every irreducible representation $R_i$ of the finite group. 
In the case $n \leq 3$, these bundles span the $K$-theory of $X$, and via the Chern
character isomorphism, $\{\ch(\sR_0), \ch(\sR_1), \ldots,
\ch(\sR_r)\}$ span $H^{*} (X; \IQ)$.

The idea is to use the Atiyah-Patodi-Singer (APS) index theorem for studying the ring
structure in cohomology. We define a generalization of the Cartan
matrix of the case $n=2$, and show a that a generalization of
Kronheimer and Nakajima's formula~\eqref{eq:dim-2-multiplication} holds:
\begin{equation}\label{eq:dim-3-multiplication}
  \left[\int_X \left(\ch(\sR_i)   - \rk(\sR_i)\right) 
               \left(\ch(R_j^\ast) - \rk (\sR_i)\right) 
  \right]_{i, j = 1, \ldots, r} = \, C^{-1}.
\end{equation}

\subsection{Overview}
In section~\ref{sec:geometry} we present a construction of Nakamura's $G$-Hilbert
scheme as a symplectic reduction of a bigger space related to the
representation theory of $G$.  Then we restrict ourselves to the case
$n=3$ and assume that the singularity $\IC^3/G$ is isolated. We prove
that the induced metric is ALE of order $6$. The decay of the metric 
allows us to use Joyce's proof of the Calabi Conjecture on ALE
manifolds to find a unique Ricci-flat ALE metric in its K\"ahler class.

In section~\ref{sec:analysis} we consider the Dirac operator and analyze its index
when completed in a weighted Sobolev norm. We compute the index of
this operator using Atiyah-Patodi-Singer index theorem.
In section~\ref{sec:vanishing} we prove a vanishing result for the index of the Dirac
operator. This result allows us to derive our geometrical
interpretation of McKay Correspondence~\eqref{eq:dim-3-multiplication}.
We conclude with remarks about future work and  related interpretations of the McKay
Correspondence.

\subsection{Acknowledgments}
This work is a continuation of the author's PhD thesis under the
supervision of Tom Mrowka. His guidance has been absolutely invaluable
for this work. We thank Mark Stern for useful conversations about 
index theory on non-compact manifolds. We thank Andrei Caldararu and
Alastair King for helping us understand the relation between the results here and
Bridgeland-King-Reid's paper. 

\section{Nakamura's $G$-Hilbert scheme as an ALE space}\label{sec:geometry}

For $G$ a nontrivial finite subgroup of $SL(n, \IC)$ we consider 
{\em Nakamura's $G$-Hilbert scheme} $\Hilb^G (\IC^n)$. This is
the scheme parameterizing $0$-dimensional subschemes $Z$ of $\IC^n$ 
satisfying the following three conditions:
\begin{enumerate}
 \item the length of $Z$ is equal to the order of $G$, $|G|$;
 \item $Z$ is invariant under the $G$-action;
 \item $H^0 (\sO_Z)$ is the regular representation of $G$.
\end{enumerate}
Let $X$ be the component of $\Hilb^G(X)$ which contains the set of
$G$-orbits $G\cdot (\IC^n \setminus \{0\})$. 
(Usually $\Hilb^G (\IC^n)$ is the union of components of fixed points 
under the $G$-action in the Hilbert scheme of $|G|$-points in $\IC^n$.)
It comes with a natural morphism
$$
  \pi: X \to \IC^n/G,
$$
the Hilbert-Chow morphism.
The variety $X$ can be described as a symplectic/GIT 
quotient~\cite{ito-nakajima:00, sar-inf:par-res, nak:lec-hil,
  kron:ale}. We review these descriptions now.

Let $G$ be a finite subgroup of $SL(n, \IC)$ and let $Q$ be its canonical $n$-dimensional
representation. Let $R$ be the regular representation of $G$: 
$$
  R = \bigoplus_{i=0}^r \IC^{n_i} \tensor R_i,
$$
where $\{R_0, R_1, \ldots, R_r\}$ is the set of irreducible representations
of $G$, with $R_0$ being the one-dimensional trivial
representation. The number $n_i$ denotes the dimension of $R_i$.
We consider the complex vector space
$$
  \sP = Q \tensor \End (R)  \oplus \Hom (R_0, R) \oplus \Hom (R, R_0).
$$
We choose an orthonormal basis of $Q$. With this choice, an element of
$Q \tensor \End (R)$ is represented as an $n$-tuple of endomorphisms
$B = (B_1, \ldots, B_n)$.
Therefore the  group $G$ acts naturally on $\sP$ via 
$$
  g \cdot ( B, i, j)  
  = 
 (g^{-1} B g,\, g^{-1} i,\, j g). 
$$
We take the $G$-invariant part of $\sP$ and denote it by $\sM$.
Inside $\sM$ we consider the variety
$$
  \sN = \{ (B, i, j) \in \sM| B \wedge B + i j = 0\}, 
$$
with $B \wedge B = \displaystyle{\sum_{\alpha < \beta} [B_{\alpha}, B_{\beta}]}$.
This is an algebraic variety (obtained from cutting the complex
vector space $\sM$ by quadratic equations); it is actually a cone in $\sM$. 

Let $GL(R)$ be the linear automorphism of $R$, and let $\sF^c\subset GL(R)$ be the  
subgroup consisting of those elements which commute with the action of
$G$ on $R$. We take the maximal compact subgroup of $\sF^c$ and denote
it by $\sF$; it consists of those elements of $U(R)$ -- unitary automorphisms of
$R$ -- which commute with the action of $G$:
\begin{equation}\label{eq:group-sf}
 \sF = U^{G} (R_0) \times U^{G} (R_1) \times \ldots \times U^G(R_{r}).
\end{equation}
% Here $R_0, R_1, \ldots, R_{r}$ denote all the irreducible
% representations of $G$.
The natural action of $\sF^c$ on $\sP$ is given by
$$
  (B, i, j) 
  \mapsto 
  (f B f^{-1}, f^{-1} i, j f), 
  \quad f \in \sF^{c}.
$$
It preserves $\sM$ and the variety $\sN$, as well as the
induced K\"ahler forms and complex structures.

For a given character $\chi: \sF \to \IC^*$ there are two quotients we
can associate to $\sN$: the GIT quotient and the symplectic quotient.
For the GIT quotient we take the complexification of the character
$\chi$ -- also denoted by $\chi$ -- and construct it as
\begin{equation}\label{eq:git-q}
  \git \sN \chi\sF^c.
\end{equation}
For the symplectic quotient we need to consider the moment map 
$\mu: \sP \to \ff^*$ determined by the action of $\sF$ on $\sP$. 
Here $\ff^*$ denotes the Lie algebra of $\sF$. Concretely, this moment map is
$$
  \mu (B, i, j) = \frac{\i}{2} \left(  [B, B^*] + i i^* - j^* j \right),
$$
where $[B, B^*] = \displaystyle{\sum_{\alpha =1}^n} [B_{\alpha}, B_{\alpha}^*]$.
The restriction of $\mu$ to $\sM$ and respectively to $\sN$ gives
moment maps for the action of $\sF$ on each of these spaces.  
If $d \chi \in \mbox{Center}(\ff^*)$, the center of the dual Lie algebra to
$\sF$, we construct the symplectic quotient as
\begin{equation}\label{eq:sym-q}
  \frac{\sN \cap \mu^{-1} ( \sqrt{-1} d \chi)}{\sF}.
\end{equation}
By a result of Kempf and Ness~\cite{kem-nes:git-sym} there is a
bijection  between the spaces obtained via the two constructions:
\begin{equation*}
  \git \sN\chi\sF^c \, \iso \,  \frac{\sN \cap \mu^{-1} ( \i d \chi)}{\sF}.
\end{equation*}
The following theorem gives the description of these quotients in two
special circumstances: when $\chi$ is the trivial character and when
$\chi$ is the determinant character.
\begin{theorem}\label{thm:orb-hilb}
 \begin{enumerate} 
   \item[(i)] Assume that the finite group $G$ acts
     freely on $\IC^n\setminus \{0\}$. Then for $\chi: \sF \to \IC^*$
     the trivial character, 
     \begin{equation}\label{eq:git-sym-orb} 
       \git\sN{}\sF^c \, \iso \, \frac{\sN \cap \mu^{-1} (0)}{\sF} \,
       \iso \, \IC^n/G.
     \end{equation} 
     as algebraic varieties.  
   \item[(ii)] For $\chi =\det: \sF \to \IC^*$ the determinant 
    character, there exists a
    bijection between the sets 
    \begin{equation}\label{eq:git-sym-cre}
      \git \sN{det}\sF^c \, \iso \, \frac{\sN \cap
        \mu^{-1}(\zeta_{\mbox{Hilb}})}{\sF} \, \iso \,  X, 
    \end{equation} 
    where $\zeta_{\mbox{Hilb}} = d \chi \in \ff^*$ and $X$ is the largest
    connected component of Nakamura's $G$-Hilbert scheme.
   \end{enumerate}
\end{theorem}
Part (i) of the Theorem was proved by Kronheimer~\cite{kron:ale}
for $n=2$, and by Sardo-Infirri~\cite{sar-inf:par-res} for
$n \geq 3$. In Sardo-Infirri's proof it is noted that 
in the case when $G$ does not act freely on $\IC^n \setminus \{0\}$, 
we just have an inclusion of $\IC^n/G$ into the GIT quotient; 
in fact the GIT quotient might have many other components.

Part (ii) of the Theorem is proved by Nakajima~\cite{nak:lec-hil}
for $n=2$, and by Ito and Nakajima~\cite{ito-nakajima:00} for $n \geq 3$.

The relation  between the above construction and crepant resolutions of
Calabi-Yau orbifold singularities is established by the following theorem:
\begin{theorem}\label{thm:hil-cre}
 Assume $n \leq 3$. Then $X$ is nonsingular, and the Hilbert-Chow morphism
 $$ 
   \pi: X \to \IC^n/G
 $$ 
 is a crepant resolution.
\end{theorem}
In the case $n=2$ (Kleinian singularities) this theorem was proved by
%Ginzburg and Kapranov~\cite{gin-kap:hil-sch} and also by
Ito and Nakamura~\cite{ito-nakamura:96}. In the case $n=3$ the theorem was
proved by Bridgeland, King and Reid~\cite{bkr} using derived
categories techniques.
\begin{corollary}[Ito and Nakajima~\cite{ito-nakajima:00}, Corollary 4.6.]\label{cor:act-free}
 Assume $n \leq 3$. Then $\sN \cap \mu^{-1}(\zetahilb)$ is
 nonsingular and $\sF$ acts freely on it, making the quotient map
 $$
   \sN \cap \mu^{-1}(\zeta_{\mbox{Hilb}}) \to X,
 $$
 into a $\sF$-principal bundle. Moreover, the bijection
 in~\eqref{eq:git-sym-cre} is an isomorphism.
\end{corollary}
\noindent
This corollary allows us to use properties of the symplectic quotients to
obtain information about the geometry of the crepant resolution $X$.

\subsection{The geometry of the Crepant resolution}
In the case $n=2$ the geometry of the crepant resolution $X$ was
completely described by Kronheimer~\cite{kron:ale} using the
hyper-K\"ahler properties of this space.
In what follows we restrict to the case $n=3$ and try to give a similar
description of the geometry of $X$, building up on previous work of
Sardo-Infirri~\cite{sar-inf:par-res}.
\begin{assumption}
We make the essential assumption that the finite group $G$ acts
freely on $\IC^3 \setminus \{ 0\}$. 
\end{assumption} 
As a consequence, $G$ is an abelian group; it implies that the Lie
group $\sF$ is also abelian and that $\mbox{Center}(\ff^*) = \ff^*$. 
Therefore we can perform the symplectic quotient 
construction~\eqref{eq:sym-q} at any point $\zeta \in \ff^*$.
\begin{remark}
 The symplectic quotient for $\zeta = 0$ is $\IC^3/G$ endowed with the
 induced orbifold metric, symplectic form and complex structure.
\end{remark}
\begin{notation}
 To make the writing easier we denote the space $\sN \cap \mu^{-1}
 (\zeta)$ by $Y_{\zeta}$ for any $\zeta \in \ff^*$. The corresponding
 symplectic quotient is $X_{\zeta} = Y_{\zeta}/\sF$.
 In particular, $Y_{\zetahilb} = \sN \cap \mu^{-1} (\zetahilb)$ and 
 $X_{\zetahilb} = X$.
\end{notation}
Corollary~\ref{cor:act-free} tells us that $Y_{\zetahilb}$ is smooth
and that $Y_{\zetahilb} \to X_{\zetahilb}$ is an $\sF$-principal
bundle. Therefore the same statement is true for $\zeta$ in a small
convex neighborhood around $\zetahilb$. We take the cone over this
convex neighborhood, and denote it by $\sC$. At each point $\zeta \in
\sC \setminus \{0\}$ the symplectic quotient $X_{\zeta} =
Y_{\zeta}/\sF$ is defined. This is because the action by dilatation 
of the positive scalars on $\sN$ induces a map $\mu^{-1} (\zeta) \to
\mu^{-1} (t^2 \zeta)$ (here we use the property of the moment map of
being quadratic on $\sP$ and therefore on $\sN$).
The quotients $X_{\zeta}$  for $\zeta \in \sC \setminus \{ 0\}$ are
smooth and diffeomorphic to $X$. However the induced symplectic
form $\omega_{\zeta}$ varies with $\zeta$. This variation was studied
by Duistermaat-Heckman \cite{dui-hec:var-coh}:
\begin{lemma}\label{lem:dui-hec}
 For each $\zeta \in \sC \setminus \{0\}$, the $\sF$-principal bundle 
 $Y_{\zeta} \to X_{\zeta}$ has a natural connection $A_{\zeta}$
 induced by the symplectic form $\omega$ on $\sN$.

 Moreover, if $R_{\zeta}$ is the curvature $2$-form associated to the
 connection $A_{\zeta}$, the variation of the symplectic form
 $\omega_{\zeta}$ is given by
 \begin{equation}\label{eq:dui-hec}
  \partial_{\lambda} \omega_{\zeta} = R_{\zeta} \lambda,
 \end{equation}
 where for $\lambda \in \ff^*$, $\partial_{\lambda} \omega_{\zeta}$
 denotes the differentiation in the direction of the constant vector
 field $\lambda$.
\end{lemma}
\begin{proof}
 Let $p_{\zeta}: Y_{\zeta} \to X_{\zeta}$ denote the quotient
 projection and $i_{\zeta}: Y_{\zeta} \to \sN$ denote the inclusion
 map. The symplectic form on $X_{\zeta}$ satisfies 
 $ p_{\zeta}^{*} \omega_{\zeta} = i_{\zeta}^{*} \omega $. 
 
 For $y \in Y_{\zeta}$ the tangent space to $Y_{\zeta}$ has the
 following decomposition:
 \begin{equation}
  T_{y} Y_{\zeta} = T_{y} (\Fy) \oplus H_y,
 \end{equation}
 where $\Fy$ is the $\sF$-orbit through $y$, and $H_y = T_{y}
 (\Fy)^{\omega} \cap T_{y} Y_{\zeta}$  with $T_{y} (\Fy)^{\omega}$
 being the symplectic complement of $T_{y}(\Fy)$ in $T_y \sN$.
 The spaces $H_y$ for $y \in Y_{\zeta}$ determine a connection on the
 principal bundle $Y_{\zeta} \to X_{\zeta}$. We want to write down the
 corresponding connection $1$-form.

 The group $\sF$ acts on $\sN$. For $\lambda \in \ff^*$ 
 let $V_{\lambda}$ denote the corresponding vector field induced on 
 $\sN$. We define a $\ff$-valued one-form $A_{\zeta}$ on $Y_{\zeta}$ 
 as follows
 \begin{equation}\label{eq:con-zeta}
  A_{\zeta} (\lambda) = - i_{\zeta}^* \left(\iota(V_{\lambda}) \omega\right).
 \end{equation} 
 We need to verify that this is the connection form corresponding to the
 horizontal subspaces $H_y$.
 
 Let $\xi \in \ff$ and let $\xi^* \in \ff$ be the corresponding
 vector field on $Y_{\zeta}$. We need that 
 $A_{\zeta} (\xi^*) = \xi$. This follows since for each $\lambda \in \ff^*$ we have
 \begin{equation}
  \<A_{\zeta} (\xi^*), \lambda\>
    = - i_{\zeta}^* \left( \iota(V_{\lambda})\omega \right) (\xi^*)
    = - \omega (\xi^*, V_{\lambda})
    =   \< d\mu (V_{\lambda}), \xi \>
    =   \< \lambda, \xi \>.
 \end{equation}
 It remains to check that $A_{\zeta} (v) = 0$ for any $ v \in
 H_y$. Since $v$ lies in the symplectic complement we have 
 $A_{\zeta} (\lambda) (v) = - \omega (V_{\lambda}, v) = 0$. Therefore
 $A_{\zeta} (v) = 0$ and $A_{\zeta}$ is the connection $1$-form
 corresponding to the horizontal distribution given by $H_y$.

 The last thing which remains to be proved is the variational formula.
 We have
 \begin{equation*}
  \begin{split}
  p_{\zeta}^{*} (\partial_{\lambda} \omega_{\zeta}) 
    & = \partial_{\lambda} (p_{\zeta}^* \omega_{\zeta})
      = \partial_{\lambda} (i_{\zeta}^* \omega)
      = i_{\zeta}^* (\sL_{V_{\lambda}} \omega) =\\
    & = i_{\zeta}^* \left(
                       d\iota(V_{\lambda}) + \iota d \omega (V_{\lambda})
                    \right) 
      = i_{\zeta}^* \left(
                      d \iota (V_{\lambda}) \omega
                    \right) = \\
    &  = d i_{\zeta}^* \left( \iota(V_{\lambda} \omega) \right)
      = d A_{\zeta} (\lambda) = \\
    &  = p_{\zeta}^* R_{\zeta} (\lambda)
  \end{split}
 \end{equation*} 
 For the last equality we use the fact that the structure group $\sF$ of the
 bundle $Y_{\zeta} \to X_{\zeta}$ is abelian. 
 Since $p_{\zeta}$ is a Riemannian submersion it follows that
 \begin{equation*}
  \partial_{\lambda} \omega_{\zeta} = R_{\zeta}(\lambda).
 \end{equation*}
\end{proof}
The complex structure on $\sN$ induces a natural complex structure on
$X_{\zeta}$ which is compatible to the symplectic structure. 
We also obtain a Riemannian metric $g_{\zeta}$ on $X_{\zeta}$
which is compatible to the complex structure and the symplectic form.
\begin{lemma}
 For $\zeta \in \sC \setminus \{0\}$, the symplectic quotient
 $X_{\zeta}$ inherits a complex structure $J$ from $\sN$ and a metric 
 $g_{\zeta}$ which are compatible with the symplectic form $\omega_{\zeta}$.
\end{lemma}
\begin{proof}
 We need to show that $H_y$ is also the orthogonal complement of
 $T_y(\Fy)$ in $T_y (Y_{\zeta})$ and that it is left invariant by 
 the action of the complex structure $J$ on $\sN$.
 
 From the proof of the previous Lemma, we have 
 $H_y = T_{y}(\Fy)^{\omega} \cap T_{y}Y_{\zeta}$.
 Let $H'_y$ denote the orthogonal complement of $T_{y} (\Fy)$ in 
 $T_{y} Y_{\zeta}$.
 If $v \in T_{y} (Y_{\zeta}$ and $\xi \in \ff$ we have
 $$
   g(Jv, \xi^*) = \omega (Jv,J \xi^*) 
                = \omega (v, \xi^*)
                = \< d\mu(v), \xi\>
                = 0,
 $$
 since $d \mu (v) = 0$.
 Therefore $\xi^* \perp J v$ and $ J \xi^* \in T_{y} Y_{\zeta}^{\perp}$.
 
 We want to check that $J$ leaves $H'_y$ invariant. Let $v \in H'_y$. 
 We have
 $$
   \< d\mu(Jv), \xi \> = \omega(Jv,\xi^*)
                       = - g (Jv, I\xi^*)
                       = - g(v, \xi^*)
                       = 0,
 $$ 
 and therefore $J v \in \Ker (d \mu)$. Thus $J v \in T_{y}
 (Y_{\zeta})$; since $J v \perp \xi^*$ it follows that $ J v \in H'_y$.
 
 From the above considerations it follows that 
 \begin{equation}\label{eq:tan_decomp}
   T_{y} (\sN) = H'_y \oplus T_{y} (\Fy) \oplus J T_{y}(\Fy).
 \end{equation}
 Since the symplectic form $\omega$ on $\sN$ is compatible with the
 complex structure it follows that $H'_y = H_y$.
\end{proof}
The induced metric $g_{\zeta}$ on $X_{\zeta}$ has a special property:
it is an ``ALE metric''.

\subsection{ALE metrics}

\begin{definition}\label{def:ale}
 A Riemannian manifold $(X, g)$ of real dimension $m$ is called 
 {\em asymptotically locally  euclidean}  (ALE)
 of order $\mu>0$ if there exists a compact set 
 $K \subset X$ and 
 a finite subgroup $G$ of $SO(m)$ so that $X \setminus K$ -- the 
 ``end'' of $X$ -- is diffeomorphic to $(\IR^m \setminus B_R)/G$
 for some $R >0$.
 Under this diffeomorphism the metric on $X \setminus K$ is of the
 form
 $$
   g_{ij} = \delta_{ij} + O(r^{-\mu}), 
   \quad 
   \nabla^k g_{ij} = O(r^{-\mu -k})  \mbox{for $k \leq 1$}, 
 $$
 in the geodesic polar coordinates $\{r, \Theta\}$
 -- the {\em ALE coordinates} -- induced on $X \setminus K$.
\end{definition}
This definition apparently depends on the choice of ALE coordinates. 
However, it can be shown~\cite{par-lee} that the ALE structure is
determined by the metric alone.
\begin{remark}
 The ALE condition translates into the fact
 that the group $G$ acts freely on $S^{m-1}$.
 We think of an ALE manifold as a manifold with boundary $Y =
 S^{m-1}/G$ at infinity. 
\end{remark}
\begin{example}
 A crepant resolution of singularities of a Calabi-Yau orbifold admits
 ALE metrics, if the dimension of the orbifold is $\leq 3$.
\end{example}  
Sardo-Infirri showed that a crepant resolution of the isolated
singularity $\IC^3/G$ admits an ALE metric. His ALE metric has order
$\mu = 4$. For our purposes (see Theorem~\ref{thm:joyce}) we need a
stronger decay on the metric: ALE of order
$\mu = 6$. We work with Sardo-Infirri's approach and show that his
estimate can be improved.
\begin{proposition}\label{prop:ale-met}
 For $\zeta \in \IC \setminus \{ 0 \}$ the metric $g_{\zeta}$
 on $X_{\zeta}$  induced via the symplectic quotient construction is
 ALE of order $\mu = 6$.
\end{proposition}
\begin{proof}
 Let $\zeta \in \sC \setminus \{0\}$.
 Via the map 
 $$
  \IR^{6} \setminus \{ 0 \} \to \IR^6/G \to X_{\zeta},
 $$
 we transfer all the structure (symplectic form, Riemannian metric,
 complex structure) on $X_{\zeta}$ to $\IR^6 \setminus \{0\}$.
 Therefore we have 
 \begin{itemize}
  \item $\{ g_{\zeta} \}_{\zeta \in \sC}$ family of metrics
    on $\IR^6 \setminus \{ 0 \}$, with $g_0 = \delta$
    the Euclidean metric.
  \item $\{ \omega_{\zeta} \}_{\zeta \in  \sC}$ family of symplectic
    forms on $\IR^6 \setminus \{ 0 \}$, with $\omega_0$ being the
    standard symplectic form on $\IC^3 \setminus \{0\}$.
  \item $\{ A_{\zeta} \}_{\zeta \in \sC}$ family of connections on the
    principal $\sF$-bundle $ ( \IR^n \setminus \{ 0 \}) \times \sF) \to
    \IR^n\setminus \{ 0\}$, with $A_0$ the trivial connection.
 \end{itemize}
  Let $(r, \Theta)$ be the polar coordinates on $\IR^6 \setminus \{ 0 \}$.
 The dependence of the family $\{ g_{\zeta}\}$ on $\zeta$ is analytic,
 and we have the following power series expansion around $0$:
 \begin{equation}\label{eq:met-exp-zeta}
   g_{\zeta}|_{r =1} 
   = 
   \sum_{|\nu| \geq 0} f_{\nu} \zeta^{\nu}.
 \end{equation}
 The rescaling via the map
 $\mu^{-1} (\zeta) \to  \mu^{-1} (t^2 \zeta)$ gives 
 $$
   g_{\zeta} (r, \Theta) = g_{r^{-2} \zeta} (1, \Theta).
 $$
 Then
 \begin{eqnarray*}\label{eq:met-exp-r}
  &g_{\zeta}(r,\Theta)&=g_{r^{-2}\zeta}(1, \Theta)
                       =\sum_{|\nu| \geq 0} 
                         f_{\nu} \zeta^{\nu} r^{-2|\nu|}\\
  &                   &=\sum_{k \geq 0} h_k(\Theta) r^{-2 k},
 \end{eqnarray*}
 where $h_k (\Theta) = \displaystyle{\sum_{|\nu| = k}} f_{\nu} \zeta^{\nu}$.

 In order to show that $g_{\zeta}$ is ALE we need that $h_0= \delta$
 and  $h_1 = h_2 =0 $.
 For $h_0$ we have
 \begin{equation*}
   h_{0} (\theta) = g_{0} (1, \theta) = \delta (1, \theta).
 \end{equation*}
 For $h_1 = 0$ we need to show that $\partial_V g_{\zeta}|_{\zeta =0} 
 (1, \theta) =0$, or  equivalently that $\partial_V
 \omega_{\zeta}|_{\zeta =0} (1, \theta) =0$. 
 Duistermaat-Heckman's formula~\eqref{eq:dui-hec} gives that
 \begin{equation}
   \partial_V \omega_{\zeta} (1, \theta) = \<V, R_{\zeta} (1, \theta)\>,
 \end{equation}
 where $R_{\zeta}$ is the curvature of the connection form $A_{\zeta}$.
 Since $R_0 = 0$ we obtain that $h_1 (\theta) = 0$.

 We have left to show that $h_2 = 0$. By \eqref{eq:dui-hec} this statement is equivalent to 
 $\partial_{\lambda} R_{\zeta}|_{\zeta = 0} = 0$.
 We have
 \begin{equation*}
  \begin{split}
    \partial_{\lambda} R_{\zeta}|_{\zeta = 0} 
     & = \frac{d}{d t}|_{t =0}  R_{\lambda t} \\
     & = \frac{d}{d t}|_{t =0}  (\sqrt t)^3 R_{\lambda} \\
     & = 0,
   \end{split}
 \end{equation*}
 since we have the rescaling  $\mu^{-1} (\lambda t) \to \mu^{-1}
 (\lambda \sqrt{t})$ and since the ambient space $X_{\zeta}$ has
 complex dimension $3$.
\end{proof}
\subsection{Tautological Vector Bundles on the Crepant Resolution}
From the definition of $X_{\zetahilb}$ there exists a natural vector bundle
\begin{eqnarray*}
  & \sR \to X_{\zetahilb} \\
  & \sR = Y_{\zeta} \times_{\sF} R, 
\end{eqnarray*}
with $R$ the regular representation of $G$. Let $\{R_0, R_1, \ldots,R_r\}$ 
be the irreducible representations of $G$ with $R_0$ the trivial
representation. The regular representation decomposes as a 
$G \times G$ module as
$$
  R = \bigoplus_i W_i \tensor R_i,
$$
where $\dim W_i = \dim R_i$.
\begin{remark}
 Since we are in the case of isolated singularities, all the $W_i$ are
    one-dimensional complex vector spaces.
\end{remark}
The vector bundle $\sR$ has a corresponding decomposition as a
$G$-module:
$$
  \sR = \bigoplus_{i} \sR_i \tensor \underline{R_i}, 
$$
where we define the vector bundles $\sR_i \to X_{\zetahilb}$  by
$$
  \sR_i = Y_{\zeta} \times_{\sF} W_i,
$$
and $\underline{R_i}$ denotes the trivial vector bundle with fiber
$R_i$.
We refer to the bundles $\sR_i$ as {\em tautological vector
  bundles} on X.
\begin{proposition} [Bridgeland-King-Reid~\cite{bkr}]\label{prop:bkr}
 The set of Chern characters 
\newline $\{ \ch(\sR_0), \ch(\sR_1), \ldots, \ch(\sR_r)\}$ forms a
 basis of $H^{*} (X, \IQ)$.
\end{proposition}
We want to figure out the ring structure in cohomology of $H^{*} (X,
\IQ)$. For this we need first to introduce an analytical set-up for
doing analysis on the crepant resolution $X$.

First we need to notice that the bundles $\sR_i$ are endowed with
natural ``asymptotically flat'' connections. We call a connection $A$ {\em 
asymptotically flat} if there exists a flat connection $A$ defined on
the ALE end $X \setminus K$ such that, under a suitable
trivialization, the two connections satisfy
$$
  A-A_0 = O(r^{-1}), \quad \nabla^k A - \nabla^k A_0 = O(r^{-1-k}).
$$
\begin{lemma}
 The connection $A_{\zeta}$ on $X_{\zeta}$ is asymptotically flat
 and the corresponding curvature behaves like $\sO(r^{-4})$ at infinity.
\end{lemma}
This was proved by Gocho and Nakajima~\cite{goc-nak} in the case $n=2$. Their proof
generalizes straightforwardly to our case. For the sake of completion we 
include it here.
\begin{proof}
 The curvature form of the connection $A_{\zeta}$ is given by
 $$
   R_{\zeta} (V, W) = - A_{\zeta} ([\widetilde{V}, \widetilde{W}]^v),
 $$
 where $\widetilde{V}$ denotes the horizontal lift of the tangent vector
 $V$ to $X_{\zeta}$. The vertical component of a tangent vector $U$ on
 $Y_{\zeta}$ is denoted by $U^{v}$.
 The bundle map $Y_{\zeta} \to X_{\zeta}$ is a Riemannian submersion.
 Therefore the Levi-Civita connections on the two spaces are related by
 $$
   \widetilde{\nabla^{^{X_{\zeta}}}_V W} =
   \nabla^{^{Y_{\zeta}}}_{\widetilde{V}} \widetilde{W}
   - \frac{1}{2} ([\widetilde{V}, \widetilde{W}]^v).
 $$
 The space $Y_{\zeta}$ is a smooth submanifold of $\sN$ and let $\Pi$
 denote the second fundamental form. Then the above formula becomes
 $$
  \widetilde{\nabla^{^{X_{\zeta}}}_V W}
  = \nabla^{^{\sN}}_{\widetilde{V}} \widetilde{W} 
   -\Pi(\widetilde{V}, \widetilde{W}) 
   - \frac{1}{2} ([\widetilde{V}, \widetilde{W}]^v).
 $$
  The previous formula applied to $I\widetilde{W}$ gives 
 $$
  \nabla^{^{\sN}}_{\widetilde{V}} I \widetilde{W}
  = \widetilde{\nabla^{^{X_{\zeta}}}_V IW}
    +\Pi(\widetilde{V}, I \widetilde{W}) 
    +\frac{1}{2} ([\widetilde{V}, I \widetilde{W}]^v).
 $$
 Since the complex structure on $\sN$ is parallel with respect to
 $\nabla^{\sN}$ and given the decomposition~\eqref{eq:tan_decomp}
 of the tangent to $T_{y} \sN$  we have
 $$
   \nabla^{^{X_{\zeta}}}_V IW = I \nabla^{^{X_{\zeta}}}_V W, \quad
   I\, \Pi(\widetilde{V}, \widetilde{W}) = \frac{1}{2} [\widetilde{V},
   \widetilde{W}]^v, \quad
   \Pi(\widetilde{V}, I \, \widetilde{W}) = \frac{1}{2} I\,
   [\widetilde{V}, I \widetilde{W}]^v.
 $$
 The first relation shows that the complex structure on $X_{\zeta}$ is
 parallel to the Levi-Civita connection. From the second identity we
 deduce that
 $$
   \Pi(\widetilde{V}, \widetilde{W}) = - \frac{1}{2} [\widetilde{V},
   \widetilde{W}]^v.
 $$
 Therefore 
 $$
   \widetilde{\nabla^{^{X_{\zeta}}}_V W} 
  = \nabla^{^{\sN}}_{\widetilde{V}} \widetilde{W}
  + \frac{1}{2} I\, [\widetilde{V}, I\, \widetilde{W}]^v
  - \frac{1}{2} [\widetilde{V}, \widetilde{W}]^v.
 $$ 
 For $\xi^*$ a vertical tangent vector to $Y_{\zeta}$ we have
 \begin{equation*}
  \begin{split}
   g([I\,\widetilde{V}, \widetilde{V}]^v, \xi^*) 
    & = \, g (I\, [I\, \widetilde{V}, \widetilde{W}]^v, I\, \xi^*)\\
    & = \, g (\nabla^{^{\sN}}_{\widetilde{V}} \widetilde{W}, I \,\xi^*)\\
    & = 2 \, \< \xi, \Hess \mu(\widetilde{V}, \widetilde{W}) \>.
  \end{split}
 \end{equation*}
 The last identity follows since
 $$
   \Pi_{- \grad \mu} (\widetilde{V}, \widetilde{W})) 
  = \< \Pi(\widetilde{V}, \widetilde{W}), \grad{\mu} \> 
  = \Hess \mu (\widetilde{V}, \widetilde{W}).
 $$
 Choosing $\xi = R_{\zeta} (V, W)$ the above formula leads to
 \begin{equation*}
  \begin{split}
   |R_{\zeta} (I V, W)|^2 
    & = \< A_{\zeta} ( [I\, \widetilde{V}, \,\widetilde{W}]^v),
           A_{\zeta} ( [I\, \widetilde{V}, \, A_{\zeta} (R_{\zeta} (I \, V, W)^*))\>\\
    & \leq |A_{\zeta}|^2 \, g([I\, \widetilde{V},\, \widetilde{W}]^v,
    R_{\zeta}(I\,V, W)^*).
  \end{split}
 \end{equation*}
 Putting everything together we obtain the following inequality
 \begin{equation}
  |R_{\zeta} (I\, V, W)| \leq 2 \, |A_{\zeta}|^2  |\Hess \mu
   (\widetilde{V}, \widetilde{W})|.
 \end{equation}
 Since the connection form scales like $A_{\zeta} (r, \Theta) = r^{-1}
 \, A_{r^{-2}\zeta} (1, \Theta)$ and since $\Hess \mu$ behaves like
 $r^{-2}$ (the metric is ALE) it follows that the curvature is of the
 form $R_{\zeta} = \sO (r^{-4})$.
\end{proof}
\begin{remark}
Since the vector bundles $\sR_i$ are associated to the principal
$\sF$-bundle $Y_{\zeta} \to X_{\zeta}$ the corresponding natural connection is
asymptotically flat. 
\end{remark}

\section{Analysis on ALE manifolds}\label{sec:analysis}
Let $(X^{m}, g)$ be an ALE space of order $\mu > 0$. From the
Definition~\ref{def:ale} this describes a Riemannian manifold with one
end which at infinity resembles the quotient $\IR^m/G$ of the
Euclidean space $\IR^m$ by a finite subgroup $G$ of $SO(m)$. The
Riemannian metric $g$ is required to approximate the Euclidean metric
up to $\sO(r^{-\mu})$. From the previous section, ALE spaces arrive as 
crepant resolutions of isolated Calabi-Yau orbifolds. 

For our purposes, we consider an ALE manifold which admits a 
$\Spin$-structure. 
(If $X$ is a crepant resolution, then $c_1(X)=0$, which
implies that $w_2(X)=0$, i.e. $X$ is spin.)
The choice of a $\Spin$-structure on the ALE end, $X\setminus K$,
is equivalent to the choice of a $G$-invariant $\Spin$-structure on 
$\IR^m \setminus B_{R}$. Since $H^{1} (\IR^m \setminus B_R; \IZ_2)$ 
is trivial, the $\Spin$-structure over this space is trivial, 
but with a $G$-action induced from the natural inclusion of $G$ into $SO(m)$.
A $\Spin$-structure on $X$ is an extension of the $\Spin$-structure 
on $X \setminus K$ to the entire $X$.

Moreover, since the ALE spaces we are dealing with arrive as crepant
resolutions of Calabi-Yau orbifolds, we  restrict ourselves to K\"ahler
ALE spaces of complex dimension $n$ ($m= 2n$). The action of the finite group
$G$ must preserve the K\"ahler structure on the end $(\IC^n
\setminus B_R)/G$, and therefore we consider 
$G \subset SU(n) \subset SO(2n)$.

On this spaces we have the two spin bundles, 
$S^{+}$ and $S^{-}$, associated to the two half-spin
representations of $\Spin(2n)$. For $E$ a Hermitian complex vector
bundle on $X$ we consider the corresponding twisted Dirac operator
$$
  D^{+}: \Omega^0(X; S^{+} \tensor E) 
           \to 
           \Omega^0(X; S^{-}\tensor E).
$$
We are interested in studying Fredholm properties of this first order
elliptic differential operator.

If the appropriate analytical setting over compact manifolds
is that of Sobolev spaces, in the case of ALE manifolds we need to take 
into consideration the behavior of the functions at infinity. Therefore
the need to work with ``weighted Sobolev spaces''.
Let $(X, g)$ be an ALE manifold, with ALE coordinates $(r, \Theta)$
on $X \setminus K$. Let $\rho: X \to [1, +\infty)$ so that 
\begin{itemize}
 \item  outside the ball of radius $2R$ we have 
        $\rho (r, \Theta) = r$;
 \item  on the compact subset $K$ of $X$ we have $\rho =1$;
 \item  on the collar $B_{2R} \setminus B_R$ we have
        $1 \leq \rho \leq 2R$.
\end{itemize}
Using $\rho$ we define the {\em $\alpha$-weighted $L^2$-norm}  for a
smooth compactly supported section of the Hermitian vector bundle $E$:
$$
  ||f||^2_{_{\scriptstyle{{L^2_{\alpha}}}}} 
  = \int_X \rho^{-2 n}\, |\rho^{\alpha} f|^2 d\vol (g).
$$ 
We denote the completion of $\Omega^{0}_{\text{comp}}(X, E)$ with
respect to the $L^2_{\alpha}$-norm by $L^2_{\alpha} (X, E)$, or 
$L^2_{\alpha}(X)$ when the bundle is clear from the context. 
A section $f$ which is bounded in the $L^2_{\alpha}$-norm behaves
like $O(\rho^{-\alpha})$ on the infinite end.

We assume now that the Hermitian complex bundle $E$ is equipped 
with a connection $A$ which is asymptotically flat. Using the
connection $A$ we extend the previously defined inner
product to an $L^2_{k, \alpha}$-norm:
$$
  ||f||^2_{_{\scriptstyle{L^2_{k, \alpha}}}}
  = \sum_{j= 0}^k \,
  \int_X \rho^{- 2 n}\, |\rho^{j + \alpha} \nabla_A f|^2 d\vol (g).
$$
We denote the corresponding completions by 
$L^2_{k, \alpha} (X, E)$ or $L^2_{k, \alpha} (X)$.

The Fredholm properties of the Dirac operator extended in a weighted
Sobolev space depend on the spectrum of the restriction to the
boundary at infinity. We recall the structure of this spectrum in the
proposition below. This result appears in \cite{baer:spectrum} and \cite{thesis}.

The spin bundles $S^{\pm}$ restricted to the boundary at infinity 
$Y= S^{2n-1}/G$  can each be identified with the spin bundle $S$ of Y
associated to the spin representation of $\Spin(2n-1)$.
\begin{proposition}\label{prop:spectrum}
 For the round metric on $S^{2n-1}$, the eigenspaces of the Dirac
 operator are
 \begin{equation}
  \begin{array}{lll}
    V_{a+\frac{1}{2} \geq \frac{1}{2} \geq \ldots \geq \frac{1}{2}}
        & \text{with eigenvalue} 
            & \scriptstyle{(-1)}^{n}\,
              \frac{2n- 1 + 2\, a }{2} \\[10 pt]
      V_{-\frac{1}{2}\geq \ldots \geq -\frac{1}{2} \geq -b - \frac{1}{2}} 
        & \text{with eigenvalue}
            & \frac{2 n - 1 + 2\, b}{2}\\[10 pt]
      V_{a + \frac{1}{2} \geq \ldots \geq \frac{1}{2}
        \geq \underbrace{\scriptstyle{-\frac{1}{2} \geq \ldots \geq -\frac{1}{2}}}_{r}
        \geq -b - \frac{1}{2}}
      &\text{with eigenvalues}   
                 & \frac{(-1)^{n+r} + 2( n + a + b)}{2}\\[5 pt]
             &   & \frac{(-1)^{n+r} - 2(n + a + b)}{2},
  \end{array}
 \end{equation}
 where $a$ and $b$ range over the positive integers.
 Here we denote by $\widetilde{U(n)}$ the double cover of the unitary
 group in order to think of the sphere $S^{2n-1}$ as the homogeneous
 space $\widetilde{U(n)}/\widetilde{U(n-1)}$. The vector space
 $V_{\mu_1 \geq \mu_2 \geq \ldots \geq \mu_n}$ represents the
 irreducible representation of $\widetilde{U(n)}$ with highest weight
 $\mu_1 \, x_1 + \mu_2 \, x_2 + \ldots + \mu_n \, x_n$ (we choose the
 fundamental Weyl chamber of $\widetilde{U(n)}$ to be $x_1 \geq x_2 \geq
 \ldots \geq x_n$).
\end{proposition}
\noindent
Observe that from the above proposition the smallest positive
eigenvalue is $\lambda^{+} = \frac{2n-1}{2}$ and the largest negative
one is $\lambda^{-} = - \frac{2n-1}{2}$.
\begin{theorem}\label{thm:index-ale}
 For the weight $\alpha$ so that $0 < \alpha < 2 n - 1$
 the closure of the Dirac operator in the $L^2_{\alpha}$-norm
 \begin{equation}\label{eq:dirac-ale}
   D^{+}_{\alpha}: L^2_{1, \alpha}\, (X, S^{+} \tensor E) 
                   \to 
                   L^2_{\alpha +1}\, (X, S^{-} \tensor E)
 \end{equation}
 is Fredholm with index given by
 \begin{equation}\label{eq:index}
   \index D^{+}_{\alpha} = \int_X \ch(E) \hat{A}(p) - \frac{\eta_E}{2}.
 \end{equation}
 Here $\eta_E (s)$ is the $\eta$-function of the Dirac operator
 restricted to the boundary at infinity $Y = S^{2n-1}/G$,
 and $\hat{A} (p)$ is the Hirzebruch $\hat{A}$-polynomial applied to 
 the Pontrjagin forms $p_i$ of the ALE metric on $X$.
\end{theorem}
%%%
The proof of this theorem is inspired by techniques from gauge theory
 \cite{mor-mro-rub}; it is  a consequence of combining the Atiyah-Patodi-Singer
index theorem for manifolds with cylindrical ends  and results of
Lockhart and McOwen, \cite{loc-mco}. 
There exist some other methods to see for which weights $\alpha$ the
Dirac operator is Fredholm, \cite{nir-wal, cho-chr, bartnik}.
The approach we choose here is convenient for our purposes of expressing the
index as a topological object.

\subsection{Proof of the Theorem~\ref{thm:index-ale}}
We change the ALE metric conformally into $\wt{g} = \rho^{-2} \, g$.
It transforms the ALE manifold into a manifold with a cylindrical end.  The
two conformally equivalent metrics give rise to the same spin bundles,
with the same Hermitian metric, but with different Clifford
multiplication.  If $\gamma: TX \to \Hom (S^{+}, S^{-})$ is the
Clifford multiplication corresponding to the metric $g$, then for
$\wt{g}$ it changes into $\wt{\gamma} = \rho \gamma$.  The Dirac
operators, $D^{\pm}$ and $\widetilde{D}^{\pm}$, corresponding to the two
metrics are related via
\begin{equation}\label{dirac_ale_to_cyl}
  \widetilde{D}^{\pm} = \rho^{\frac{2 n+1}{2}} 
                    D^{\pm}
                    \rho^{ -\frac{2 n-1}{2}} \, .
\end{equation}
Under the conformal change the compact set $K$ 
does not encounter any change in its geometry, but the
ALE end $\IR^{2n} \setminus B_R$ changes into the cylindrical end 
$[T, +\infty) \times Y$, where $T = \log R$.
We take $\tau = \log \rho: X \to [0, +\infty)$. After a
reparametrization we can assume that $\tau$ has the following properties:
\begin{itemize}
 \item on $[3 T, + \infty) \times Y$,  $\tau$ agrees with the natural
       projection onto $[3 T, +\infty)$;
 \item on $[2 T, 3 T]$, $\tau$ is between $0$ and $1$;
 \item on $X \setminus [2 T, + \infty)$, $\tau \equiv 0$.
\end{itemize}
A function which behaves as $O(\rho^{-\alpha})$ at infinity 
behaves also like $O(e^{- \alpha \tau})$. Then the corresponding
weighted Sobolev norm for the metric $\widetilde{g}$ is 
$$
  ||f||^2_{_{\scriptstyle{\widetilde{L}^2_{\alpha}}}}
   = 
  \int_X |e^{\alpha \tau} f|^2 d\vol (\widetilde{g}).
$$
We denote by $\widetilde{L}^2_{\alpha}(X)$ the corresponding weighted Sobolev space.
The Dirac operator \eqref{dirac_ale_to_cyl} on the Riemannian manifold $(X, \widetilde{g})$
acts in the following manner
\begin{equation}
 \widetilde{D}^{\pm}_{\alpha - \frac{2 n-1}{2}}: 
       \widetilde{L}^2_{\alpha - \frac{2 n-1}{2}}\, (X; S^{\pm} \tensor E)
       \to
       \widetilde{L}^2_{\alpha - \frac{2 n-1}{2}}\, (X; S^{\mp} \tensor E).
\end{equation}
Let $\beta = \alpha - \frac{2 n-1}{2}$. We have replaced the problem of
studying Fredholm properties of the Dirac operator on a ALE manifold
with the similar problem on a manifold with cylindrical end.
The study of the operator~\eqref{eq:dirac-ale} is equivalent to the
study of
\begin{equation}\label{eq:dirac-cyl}
 \widetilde{D}^{\pm}_{\beta}: \widetilde{L}^2_{\beta}\,(X; S^{\pm}\tensor E)
                              \to
                              \widetilde{L}^2_{\beta}\, (X; S^{\mp}\tensor E).
\end{equation}
We can further reduce the problem just to the study of Fredholm
properties on a manifold with the product cylindrical metric on the 
infinite end. The following lemma takes care of the weighted
Sobolev spaces.
\begin{lemma}
 A metric $\widetilde{g}$ which arrives as the conformal transform of
 an ALE metric, is almost the product metric  $\widetilde{\delta} = dt^2 + d \Theta^2$
 on the cylinder.
 As a consequence the $\beta$-weighted Sobolev norms are equivalent 
 and therefore the Hilbert space completions in these norms are the
 same. 
\end{lemma}
\begin{proof}
 An orthonormal frame for the product metric on the cylinder
 is $\{dt, d\theta_i\}$.  Since
 $$
   |\widetilde{g}-\widetilde{\delta}| \leq e^{-\mu \tau},
 $$
 an orthonormal frame for $\widetilde{g}$ has the form
 $\{dt, h_{ij} d\theta_{ij}\}$ where $h_{ij} = 1 + O(e^{-(\mu + 2)
 \tau)})$, and the volume form is 
 $d \vol(\widetilde{g}) = d \vol(\widetilde{\delta}) + O(e^{-\mu \tau})$.
 Then,
 \begin{eqnarray*}
   & \int_X |e^{-\beta \tau} f|^2 d\vol (\widetilde{g}) 
      & \leq \int_X |e^{-\beta \tau} f|^2 d\vol  (\widetilde{\delta})
             +
             \int_X O(e^{-\mu \tau}) |e^{-\beta \tau} f|^2
                    d\vol (\widetilde{\delta})\\
   &  & \leq C \; \int_X |e^{-\beta \tau} f|^2 d\vol (\widetilde{\delta}).
 \end{eqnarray*}
 The other way around goes the same.
\end{proof}
In order to complete the reduction to a manifold with the product
metric on the cylindrical end we use the following result of Lockhart
and McOwen~\cite{loc-mco}.
\begin{proposition}[Lockhart and McOwen]
 Let $P_{\beta}$ be the  Dirac operator associated to the product
 metric on the cylindrical end of $X$.
 Then $P_{\beta}$ and  $\widetilde{D}_{\beta}$  are Fredholm  for exactly the
 same values $\beta$. Moreover, their Fredholm indices are equal,
 $ \index P_{\beta} = \index \widetilde{D}_{\beta}$.
\end{proposition}
\begin{proof}
 Assume that $P_{\beta}$ is Fredholm. Let $t$ be a positive real
 number. The idea is to cut and paste the operator $P_{\beta}$ to a 
 new operator $P^{'}_{\beta} = P_{\beta} + (1-\phi_t) (D_{\beta} -
 P_{\beta})$. Here $\phi_t$ is a smooth function
 which is $1$ on $X_t$ and has support in $X_{2  t}$ (notation: $X_{t} = X
 \setminus [t, +\infty) \times Y$). Since the space of Fredholm operators is
 open, for $t$ very large the operator $P^{'}_{\beta}$ is Fredholm and
 has $\index P^{'}_{\beta} = \index P_{\beta}$.
 This new operator $P^{'}_{\beta}$ has the property that for $t > 3T$ is
 $\widetilde{D}_{\beta}$.
 Performing again a cutting and pasting procedure we obtain a parametrix for
 $\widetilde{D}_{\beta}$, constructed out of a parametrix for
 $\widetilde{D}_{\beta}$ on $X_{4 t}$, for example, and a parametrix for 
 $P^{'}_{\beta}$.
 Therefore $\widetilde{D}_{\beta}$ is Fredholm.

 The operator norm between $\widetilde{D}_{\beta}$ and $P_{\beta}$ 
 is exponentially decaying. Therefore the  two operators belong to the
 same connected set in the space of Fredholm operators and thus their
 Fredholm indices are equal.
\end{proof}
 On the cylindrical end the operator $\widetilde{D}$ has the form
 $$
   \widetilde{D} = \frac{\partial}{\partial t} + B + O(e^{-\mu t}),
 $$
 where $B$ is the Dirac operator restricted to the boundary at
 infinity $Y = S^{2n-1}/G$.
 Lockhart and McOwen proved that the operator $P_{\beta}$ (and
 therefore $\widetilde{D}_{\beta}$) is Fredholm for those weights $\beta$ for
 which the operator $B - \beta\, I$ is invertible. In other words, if $\beta$ is not
 an eigenvalue for $B$ then $P_{\beta}$ is Fredholm.
 In the case when $B - \beta \, I$ is not invertible problems arise since
 this case the operator $P_{\beta}$ does not have closed range, which
 is an essential condition for the operator to be Fredholm.

In order to exhibit the index of $P_{\beta}$ we employ the
Atiyah-Patodi-Singer theorem. First we need to establish the set-up:
%%%%%
we consider the manifold with boundary $X_{2 T}$ with the metric which
is the product metric on $[T, 2T]\times Y$. The Dirac operator 
$$
  P: \Omega^0 (X_{2 T}, \,  S^{+} \tensor E) 
     \to 
     \Omega^0 (X_{2 T}, \,  S^{-} \tensor E).
$$ 
has the following form on the cylinder $[T, 2 T] \times Y$:
$$
  P = \frac{\partial}{\partial t} + B,
$$
where $B$ is the Dirac operator on the boundary $Y$. Its $L^2$-adjoint is
$$
  P^{*}= \frac{\partial}{\partial t} - B.
$$
We consider the spectral projection
$$
  \Pi^{+}: \Omega^0 (Y, \, S \tensor E) 
           \to 
           \Omega^0 (Y,\, S \tensor E)
$$
onto the span of eigenvectors corresponding to the positive eigenvalues of
$B$. The space of all smooth sections which satisfy $\Pi^{+} f(2 T, \,
\cdot \,)$ we denote by $\Omega^0 (X_{2 T}, \, S^{+}\tensor E; \,
\Pi^{+})$. Its completion in the $L^2_k$-norm we denote by
$\Omega^0_{k} (X, \, S^{+}\tensor E; \, \Pi^{+})$.
We take the closure of the Dirac operator $P$ on $L^2$ with
domain given by the global boundary condition induced by $\Pi^{+}$:
\begin{equation}\label{aps}
  \sP: \Omega^{0}_{1} (X_{2 T},\, S^{+} \tensor E; \, \Pi^{+}) 
       \to 
       \Omega^{0}_{0}(X_{2 T},\, S^{-} \tensor E).
\end{equation}
The Atiyah-Patodi-Singer index theorem states that the 
operator~\eqref{aps} is Fredholm; its index is 
$$
  \index \sP = \int_{_{\scriptstyle{X_{2T}}}} \ch(E)\, \hat{A}(p) - \frac{\eta_E}{2}.
$$
To conclude the proof of Theorem~\ref{thm:index-ale} we need the following proposition:
\begin{proposition}
 For $- \frac{2n-1}{2} < \beta < \frac{2n-1}{2}$ the operator $P_{\beta}$ is Fredholm and its
 Fredholm index is 
 $$
  \index P_{\beta} = \index \sP.
 $$
\end{proposition}
\begin{proof}
 We assume that the metric is a product metric on the cylinder 
 $[T, +\infty) \times Y$. 
 The pair of operators $P_{\beta}$ and $\sP$ can be put together into an
 excision data 
 $$
   (Z_1 = X, \quad 
    U_1 = X_{2T-\, \epsilon}, \quad
    V_1 = (T, +\infty) \times Y, \quad
    Q_1= P_{\beta}, \quad
    Q_1^{'} = \sP),
$$
where
$$
  P_{\beta}: \widetilde{L}^2_{1, \beta}\, (X) \to
  \widetilde{L}^2_{\beta}\, (X),
$$
and
$$
  \sP: L^2_{1} (X_{2 T}, \Pi^{+}) \to L^2 (X_{2 T}).
$$
Since $\tau = 0$ on  $X_{2 T}$, it follows that 
$$
  P_{\beta}|_{U_1} = \sP|_{U_1}.
$$

The second excision data we need is
$$
  (Z_2 = [T, +\infty) \times Y, \quad
   U_2 = [T, 2 T -\epsilon] \times Y, \quad
   V_2 = (T, +\infty) \times Y, \quad
   Q_2, \quad
   Q_2^{'}),
$$
where $Q_2$ is the restriction of $P_{\beta}$ to
$$
  Q_2: \left\{f\in \widetilde{L}^2_{1,\beta}(Z_2) \; |\; 
         \left( 
            1-\Pi^{+}
         \right) f(T, \, \cdot \, )=0 \right\}
       \to
       \widetilde{L}^2_{\beta}(U_2),       
$$ 
and $Q_2^{'}$ is the restriction of $\sP$ to
$$
  Q_2^{'}: \left\{f\in L^2_{1}([T, 2 T] \times Y)|
                \, \left(1- \Pi^{+}\right) f(T,\, \cdot \,) = 0, \, \,
                \Pi^{+} f (2 T,\, \cdot \,) =0 \right\}
        \to
        L^2 ([T, 2 T] \times Y).
$$
Again since $\tau$ is $0$ on $[T, 2 T] \times Y$, we have 
$$
  Q_2|_{U_2} = Q^{'}_2|_{U_2}.
$$

The two pairs of excision data are related in the following way
\begin{eqnarray*}
  & V_1 = V_2, \\
  & Q_1|_{V_1} = Q_2|_{V_2}\\
  & Q^{'}_1|_{(T, 2 T] \times Y} = Q^{'}_2|_{(T, 2 T] \times Y}.
\end{eqnarray*}
The excision principle implies that
$$
  \index Q_1 - \index Q^{'}_1 = \index Q_2 - \index Q^{'}_2,
$$
i.e.
$$
  \index P_{\beta} - \index \sP = \index Q_2 - \index Q^{'}_2,
$$
\begin{claim}
 For the operators $Q_2$ and $Q^{'}_2$ the following are true
 \begin{eqnarray*}
   & \dim \Ker Q_2 = 0,  & \quad \dim \Coker Q_2 = 0\\
   & \dim \Ker Q^{'}_2 = 0, & \quad \dim \Coker Q^{'}_2 = 0.
 \end{eqnarray*}
\end{claim}
 We sketch the proof for $\dim \Ker Q_2 =0$; the other statements
 follow in a similar way.
 Functions on $[T, +\infty) \times Y$ can be decomposed according to
 the eigenspaces of the Dirac operator $B$ on $Y$:
 $f(t, y) = \sum_{\lambda} f_{\lambda}(t) \phi_{\lambda}(y)$, where
 $\phi_{\lambda}$ are the eigenfunctions of $B$.
 Assume $f(t, y) \in \Ker Q_2$. This is equivalent to 
 \begin{equation}
  \left\{
    \begin{array}{l}
      \frac{d \, f_{\lambda}}{dt} + \lambda f_{\lambda} = 0,\\
      \mbox{with the boundary conditions:}\\
       f_{\lambda} (T) = 0, \mbox{for $\lambda < 0$}\\
       f_{\lambda} = O(e^{-\beta t}) \mbox{as $t \to +\infty$},
      \mbox{for any $\lambda$}.
    \end{array}
  \right.
 \end{equation} 
 For $\lambda < 0$ the solution to this system is 
 $$ 
   f_{\lambda} (t) = C \int_{T}^t e^{-\lambda u} du,
 $$
 which satisfies the required boundary condition if $\beta >
 -\lambda$. Since we choose $ - \frac{2 n -1}{2} < \beta$, then the previous
 inequality never happens, so the only solution is the trivial one.
 For $\lambda > 0$ the solution to the ODE is
 $$
   f_{\lambda} (t) = - \int_{t}^{+\infty} e^{-\lambda u} du
 $$
 which satisfies the decay condition only for $\beta < \frac{2n-1}{2}$.

Therefore for $-\frac{2n-1}{2} < \beta < \frac{2n+1}{2}$ the operator
$P_{\beta}$ is Fredholm. Since $\beta = \alpha - \frac{2n-1}{2}$, it
means that for $0<\alpha <2n-1$ the Dirac operator $D_{\alpha}$ in
\eqref{eq:dirac-ale} is Fredholm.
\end{proof}

\begin{remark}
 It can be proved that the Dirac operator \eqref{eq:dirac-ale} is Fredholm
 for all weights $\alpha = \beta - \frac{2n-1}{2}$ so that $\beta$ is not an
 eigenvalue of $B$ (in which case the Dirac operator does not have
 closed range). The index is going to jump -- up or down -- every time
 we cross an eigenvalue; the jump is by the dimension of the eigenspace
 corresponding to that eigenvalue.
\end{remark}

\section{A Vanishing Results on the Crepant Resolution}\label{sec:vanishing}
The index formula~\eqref{eq:dirac-ale} holds for any ALE manifold with
a metric of order $\mu \geq 1$. We are now interested in studying the
case of an ALE metric which is Ricci-flat.

In Section~\ref{sec:geometry} an ALE metric on the crepant resolution
of the isolated orbifold singularity $\IC^n/G$ was
constructed. Computations of Sardo-Infirri~\cite{sar-inf:par-res}
show that this metric is not Ricci-flat (the metric induced on
$\sO_{\IC P^2} (-3)$ the crepant resolution of $\IC^3/G$ is not
Ricci-flat; however this space has a Ricci-flat ALE metric, the
Eguchi-Hanson metric). In the situation when the ALE metric $g$ has
order $\mu = 2 n$, there is a result which allows us to overcome this problem:
Joyce's proof of Calabi's Conjecture for ALE spaces.
\begin{theorem}[Joyce~\cite{joyce:book}, Theorem $8.2.3$]\label{thm:joyce}
 Let $G$ be a finite group of $SU(n)$ acting freely on $\IC^n
 \setminus \{0\}$, and $(X, \pi)$ a crepant resolution of $\IC^n/G$. 
 Then in each K\"ahler class of ALE K\"ahler metrics of order $2 n$ on $X$ there is a
 unique Ricci-flat ALE K\"ahler metric $g$.
\end{theorem}

Assume now that $X$ is the crepant resolution of $\IC^3/G$ given by
Nakamura's $G$-Hilbert scheme. We proved in
Proposition~\ref{prop:ale-met} that $X$ comes endowed with a natural ALE
K\"ahler metric of order $6$. By
Joyce's result we have a unique Ricci-flat ALE metric in its K\"ahler
class. It is for this metric that we want to apply the index 
formula~\eqref{eq:dirac-ale}. We obtain the following vanishing result:
\begin{theorem}\label{thm:vanishing}
 Assume $G$ is a finite subgroup of $SL(3, \IC)$ which acts freely on
 $\IC^3 \setminus \{0\}$. Let $(X, \pi)$  be the crepant resolution
 given by the Nakamura's $G$-Hilbert scheme, and $g$ be the natural
 Ricci-flat ALE metric on it.
 Let $E$ be a self-dual holomorphic Hermitian bundle ($E = E^*$)
 endowed with a connection which is asymptotically flat.
 Then, for the weight $\alpha$ so that $0 < \alpha < 5$
 the Dirac operator
 $$
   D^+_{\alpha}: L^2_{1, \alpha} (X, S^+ \tensor E) 
                \to 
                 L^2_{\alpha + 1} (X, S^{-}  \tensor E),
 $$
 has vanishing index.
\end{theorem}
\begin{proof}
 The $L^2$-adjoint of the Dirac operator is
 $$
   D^{-}_{-\alpha + 5}: L^2_{1, -\alpha + 5} (X, S^{-}\tensor E)
                            \to
                           L^2_{-\alpha + 6} (X, S^{+} \tensor E)
 $$

 On a complex K\"ahler manifold, the Dirac operator and its adjoint are deeply
 related to the Dolbeault operator: we have
 $S^{+} = \Lambda^{0,0} (X) \oplus \Lambda^{0, 2}(X)$ and 
 $S^{-} = \Lambda^{0,1} (X) \oplus \Lambda^{0, 3}(X)$. The Dirac operator is 
 $$
   D_{A}^{+} = \sqrt{2} \left(
                              \delbar_{A} + \delbar_{A}^*
                        \right).
 $$
 Here $\delbar_{A}^{*}$ denotes the formal adjoint of $\delbar_A$; it is
 given by $\delbar_A^{*} = - *_{E} \delbar_A *_{E}$, where $*_E$ is the usual
 Hodge star operator associated to the metric $g$:
 $*_E: \Omega^0 (X, \,\Lambda^{p, q} \tensor E) \to \Omega^0 (X, \,
 \Lambda^{n-p, n-q} \tensor E^{*})$.
 Moreover since the metric $g$ is Ricci-flat, the bundle $\Lambda^{n,
   0}$ is trivial.  It gives a natural isomorphism $\Lambda^{n, n-k}\,
 \iso \, \Lambda^{0, n-k}$.
 Everything is encapsulated in the following diagram:
 \medskip
 \begin{equation*}
  \xymatrix{
       \Omega^{0} (X; \Lambda^{0, \mathrm{odd}} \tensor E^*)  
          \ar[r]^{*_E}_{\iso} \ar@{-->}[d]_{D^{-}}
     & \Omega^{0}(X; \Lambda^{n, n-\mathrm{odd}}\tensor E^*)
           \ar[d]^{\delbar + \delbar^*} \ar[r]^{\iso}
     & \Omega^{0} (X; \Lambda^{n-\rm{odd}} \tensor E)
           \ar[d]^{D^{+}}
            \\
       \Omega^{0} (X; \Lambda^{0, \mathrm{odd} \pm 1} \tensor E)
                 \ar[r]^{*_E}_{\iso}
     & \Omega^{0} (X;  \Lambda^{n, n- \mathrm{odd}\pm 1} \tensor E^*)
                 \ar[r]^{\iso}
     & \Omega^{0} (X; \Lambda^{n-\rm{odd} \pm 1} \tensor E)
          }
 \end{equation*}
 Completing the diagram in the $L^2_{\frac{5}{2}}$-norm,  it follows
 that
 \begin{equation}
   D^{+}_{\frac{5}{2}} = D^{-}_{\frac{5}{2}}.
 \end{equation}
 On the other hand we have that the indices of the two operators are
 related by
 $$
  \index D^{+}_{\frac{5}{2}} = - \index D^{-}_{\frac{5}{2}} 
 $$
 since  $D^{-}_{\frac{5}{2}}$ is the $L^2$-dual of
 $D^{+}_{\frac{5}{2}}$.
 Therefore we have the desired vanishing result
 \begin{equation}
  \index D^{+}_{\frac{5}{2}} = 0.
 \end{equation}
 Since the index of the Dirac operator $D_{\alpha}$ is
 constant for the weight between two critical values, it follows
 that the vanishing holds for all the weights between $0$ and $5$.
\end{proof}
\begin{remarks}
 \begin{enumerate}
   \item The proof is valid for all other crepant resolutions of
         $\IC^3/G$. They are obtained from $X$ via a flop. The flop changes the
         geometry of the exceptional fiber but since the singularity is
         isolated it does not affect the geometry at infinity.
   \item Our proof also works for $n$ any odd number, provided that we know that
         a crepant resolution exists and it has an ALE metric of order
         $2 n$.
 \end{enumerate}
\end{remarks}

\section{Geometrical McKay Correspondence}\label{sec:geom-mckay}
We want to apply the index formula~\eqref{eq:dirac-ale} together with
the vanishing result of the Theorem~\ref{thm:vanishing} to obtain a
geometrical interpretation of the McKay Correspondence similar to the
one obtained by Kronheimer and Nakajima~\cite{KN} for the case of
surface singularities. The remaining needed piece of the puzzle is
the $\eta$-invariant term.
\begin{proposition}
 Let $X$ be a crepant resolution of the isolated orbifold singularity
 $\IC^3/G$. Consider $\sE$ a Hermitian bundle over $X$ with fiber at
 infinity $E$. Then the corresponding $\eta$-invariant is 
 \begin{equation}\label{eq:eta}
  \eta_E = \frac{1}{|G|} \sum_{\substack{g \in G \\ \\ g \neq I}} \,\,
                         \frac{\chi_E (g)}
			      {- \chi_{Q} (g) + \chi_{\Lambda^2 Q}(g)}.
 \end{equation}
% \begin{equation}\label{eq:eta}
%  \eta_E = \frac{1}{|G|} \sum_{\substack{g \in G \\ \\ g \neq I}} \,\,
%                         \frac{\Trace (g|_{E})}
%			      {- \Trace(g|_{\IC^3}) + \Trace(g|_{\Lambda^2 (\IC^3)})}.
% \end{equation}
 In this formula $Q$ represents the
 $3$-dimensional representation of $G$ induced by its inclusion into 
 $SL(3, \IC)$. Also for a representation $V$ of $G$ we denote by $\chi_V$
 the character of $V$.
\end{proposition}
This proposition is a consequence of the Lefschetz fixed-point
formula, in the sense that $\eta_E$ is the contribution from the fixed
locus under the action of $G$ on $\IC^3$. It can be also proved using
the definition of the $\eta$-invariant as the analytic continuation at $0$
of the $\eta$-series corresponding to the spectrum of the Dirac
operator on the boundary at infinity of the orbifold $\IC^3/G$. 
This last approach gives the generalization of the above formula to 
the case of non-isolated singularities~\cite{thesis}.
\begin{definition}\label{def:cartan}
 Let $\{R_0, R_1, \ldots, R_r \}$ be the irreducible representations
 of the finite group $G$ (here $R_0$ denotes the $1$-dimensional 
 trivial representation).
 We consider the tensor products $Q \tensor R_i$ and $\Lambda^2 Q
 \tensor R_i$ and decompose them into irreducible representations:
 \begin{equation}\label{eq:cartan}
  \begin{split}
    & R_i \tensor \, Q = \sum_{j = 0}^r a_{ij} R_j \\
    & R_i \tensor \, \Lambda^2 Q = \sum_{j = 0}^r b_{ij} R_j.
  \end{split}
 \end{equation}
 We call the matrix 
 \begin{equation*}
   \widetilde{C}= \left[ a_{ij} - b_{ij}\right]_{i,j= 0, \ldots r},
 \end{equation*} 
 the {\em extended Cartan matrix} associated to the finite group $G$.
 By removing the first row and the first column of $\widetilde{C}$ we obtain a
 new matrix, $C$, which we call the {\em Cartan matrix} associated to
 the finite group $G$.
\end{definition}
%%%%
\begin{remark}
 This definition is a generalization of the classical 
 McKay Correspondence~\eqref{eq:mckay}. 
\end{remark}
Now we are ready to state our geometrical interpretation of the McKay
Correspondence, which is a consequence of the Vanishing Theorem~\ref{thm:vanishing}.
\begin{corollary}
 Assume that the finite subgroup $G$ of $SL(3, \IC)$ acts freely on
 $\IC^3 \setminus \{ 0\}$, and let $(X, \pi)$ be the crepant
 resolution of $\IC^3/G$ constructed using Nakamura's $G$-Hilbert scheme.
 Let $\{\sR_0, \sR_1, \ldots, \sR_r\}$ be the tautological bundles on $X$ corresponding to the
 irreducible representations of $G$.

 Then, the elements $\{\ch (\sR_1), \ldots, \ch(\sR_r) \}$ satisfy 
 the following multiplicative relations
 \begin{equation}\label{eq:mult-formula}
    \left[ 
        \int_X \left(\ch(\sR_i) - \rk(\sR_i)\right)
               \left(\ch(\sR_j^*) - \rk(\sR_j)\right)
    \right]_{i, j =1, \ldots, r}
    = \; C^{-1}.
 \end{equation}
\end{corollary}
%%%%
\begin{remark}
 It is straightforward to see that since 
 $G \subset SL(3, \IC)$, the Cartan matrix $C$ is invertible.
\end{remark}
%%%%%
\begin{proof}
 The proof of Kronheimer and Nakajima goes through without major problems:
 We consider the bundle $\sR$ which has $R_{\mbox{reg}}$, the regular
 representation of the finite group, as fiber at infinity.
 We apply our Index Theorem~\ref{thm:index-ale} to the bundle
 $\sR \tensor \sR^{*}$ and the Dirac operator completed in the
 weighted Sobolev norm with $0 < \alpha < 5$:
 \begin{equation}
  \index D^{+}_{\sR \tensor \sR^{*}} 
  = 
  \int_X \ch(\sR \tensor \sR^*) \hat{A} (p)
  -
  \frac{\eta_{R_{reg} \tensor R_{reg}^{*}}}{2}.
 \end{equation}
 Since $\sR \tensor \sR^{*}$ is a self-dual bundle, it follows
 according to the vanishing result, Theorem \ref{thm:vanishing}, 
 that the left-hand-side of the above
 formula is zero.

 The bundle $\sR$ comes with a  $G$-action on it. We decompose it
 under the $G$-action:
 \begin{equation}
   \sR = \bigoplus_{i=0}^r \sR_i \tensor \underline{R_i}.
 \end{equation}
 With respect this $G$-action the index formula becomes
 \begin{equation*}
  0 = \int_X \ch(\sR_i \tensor \sR_j^{*}) \hat{A}(p)
      - \frac{\eta_{R_i \tensor R_j^*}}{2},
 \end{equation*}
 or equivalently
 \begin{equation*}
     \int_X \ch(\sR_i \tensor \sR_j^{*}) \hat{A}(p)
      = 
     \frac{\eta_{R_i \tensor R_j^*}}{2}.
 \end{equation*}
 Multiplying on the right by the Cartan matrix $\widetilde{C} =
 \left[c_{ij}\right]_{i, j = 0, \ldots r}$, the left-hand-side becomes
 \begin{eqnarray*}
  \lefteqn{\sum_{k=0}^r c_{ik} \int_X \ch(\sR_k \tensor \sR_j^*)\hat{A}(p) =}\\
   &&  =  \sum_{k=0}^r c_{ik}\int_{X}\ch(\sR_k) \hat{A}(p)
         + \sum_{k=0}^r c_{ik}\int_{X} \ch(\sR_j^*)\hat{A}(p)
         + \sum_{k=0}^r c_{ik}\int_{X} (\ch(\sR_k) -1)(\ch(\sR_j^*) -1)\\
   && =  \sum_{k=0}^r  c_{ik} \frac{\eta_{R_k}}{2} 
         + \sum_{k=0}^r  c_{ik} \frac{\eta_{R_j^*}}{2} 
         + \sum_{k=0}^r c_{ik} \int_{X} (\ch(\sR_k) -1)(\ch(\sR_j^*) -1).
 \end{eqnarray*}
 To figure out what happens to the right-hand-side, we see
 that in terms of the corresponding characters, the relations (\ref{eq:cartan}) give
 \begin{equation*}
   \left(
     \chi_{Q} (g) - \chi_{\Lambda^2 Q} (g) 
   \right) \,
   \chi_{R_i} (g) \,
   = \, \sum_{j =0}^r c_{ij}  \, \chi_{R_j} (g).
 \end{equation*}
 We multiply by $\chi_{R_k^*} (g)$ on both sides and then we sum 
 after all $g \neq 1$ to obtain
 \begin{equation}\label{eq:eta-char}
     \sum_{g \neq 1} \chi_{R_i} (g) \chi_{R_k^*}(g)
   = \sum_{g \neq 1} c_{ij} 
                    \frac{\chi_{R_j}(g) \chi_{R_k^*}(g)}
                         {\chi_{Q}(g)-\chi_{\Lambda^2 Q}(g)} 
 \end{equation}
 Therefore the right-hand-side of the equation~\eqref{eq:eta-char} is 
 \begin{equation*}
  \sum_{j=o}^r c_{ij} \eta_{\sR_j \tensor \sR_k^{*}}
   = \frac{1}{|\Gamma|} (\delta_{ik} -1).
 \end{equation*}

 Putting these formulae together, and using again the index theorem,
 we have
 \begin{equation*}
    \frac{1}{2 |\Gamma|} (\delta_{ij}-1)
    = \frac{1}{2 |\Gamma|} (\delta_{i 0} -1) 
     + \sum_{k=0}^r     
                   c_{ik}
                   \int_{X}
                           (\ch(\sR_k) -1)
                           (\ch(\sR_j^*) -1).
 \end{equation*}
 For $i \neq 0$ it gives
 \begin{equation*}
   \sum_{k=0}^r  
               c_{ik}
               \int_{X}
                       (\ch(\sR_k) -1)
                       (\ch(\sR_j^*) -1)
   = \delta_{ij}.
 \end{equation*}
 Since $\ch(\sR_0) -1 = 0$, after inverting the above relation we obtain
 \begin{equation*}
   \int_X (\ch(\sR_k) - 1) (\ch(\sR_j^*) - 1) = \, \left(C^{-1}\right)_{kj},
 \end{equation*}
 for $k, j \neq 0$.
\end{proof}
\subsection{Comments}
 By the result of Bridgeland, King and Reid (see
 Proposition~\ref{prop:bkr}) the elements 
 $\{\ch(\sR_0), \linebreak \ch(\sR_1), \ldots, \ch(\sR_r)\}$ form a basis
 in $H^{*} (X, \IQ)$. The multiplication relations in
 $H^*(X, \IQ)$ given by~\eqref{eq:mult-formula} 
 are the generalization of the Kronheimer and Nakajima's ones.
 However, in our case they do not describe the entire ring structure in
 cohomology. For example the formula does not give any information
 about terms like $\int_X c_1 (\sR_i)^3$.   

Recent work of Craw and Ishii~\cite{craw-ishii} can be used to show
that the multiplication formula~\eqref{eq:mult-formula} holds for
any crepant resolution of an isolated singularity $\IC^3/G$. In their
work, they provide a description of any (projective) crepant
resolution of $\IC^3/G$, for $G$ abelian, as a moduli space of
$G$-constelations. Our approach from Section~\ref{sec:geometry}
applies to this moduli space: it can be described as the symplectic
reduction $X_{\zeta} = (\sN \cap \mu^{-1} (\zeta))/\sF$ at some point
$\zeta\in \ff^*$. Proposition~\ref{prop:ale-met} generalizes
straightforwardly and we obtain an ALE metric of order $\mu= 6$.
Moreover the irreducible representations of the finite group $G$ give
rise to tautological sheaves $\sR_{\zeta, i}$ which form a basis for
the $K$-theory of $X_{\zeta}$. Their Chern characters satisfy the
multiplication formula~\eqref{eq:mult-formula}.

In our view, the essence of the McKay Correspondence is: how much of
a crepant resolution of $\IC^3/G$ can be described in terms of the
finite group $G$ and its embedding in $SL(3, \IC)$? From the
discussion above all the crepant resolutions of $\IC^3/G$ have a part
of the cohomology ring which is the same -- it is given by the
multiplication formula~\eqref{eq:mult-formula}. The remaining part of
the cohomology ring should depend on the crepant resolution in the
sense that it should change under a flop. The complete description of
the cohomology ring should be carried out in future work.
\subsection{Relationship to other Results}
We present the relation of our multiplicative formula to the result of
Ito and Nakajima~\cite{ito-nakajima:00}, and Bridgeland, King and
Reid~\cite{bkr}.
The author's understanding of the following is due to a discussion
with Andrei Caldararu and Alastair King at the Isaac Newton Institute
in the summer of $2002$.

Let $\pi: X \to \IC^3/G$ be the crepant resolution of $\IC^3/G$ given
by Nakamura's $G$-Hilbert scheme.We denote by $K(X)$ the Grothendieck group of coherent
$\sO_{X}$-sheaves over $X$. The coherent sheaves which are supported
on the exceptional divisor $\pi^{-1} (0)$ generate another group which
we denote by $K_c(X)$. This can be thought of as the Grothendieck
group of bounded complexes of algebraic vector
bundles on $X$ which are exact outside $\pi^{-1} (0)$.

The result of Bridgeland, King and Reid~\cite{bkr} implies that
$\{\sR_0, \sR_1, \ldots, \sR_r\}$ form a
basis of $K(X)$. A procedure of Ito and
Nakajima~\cite{ito-nakajima:00} 
gives a basis $\{\sS_0, \sS_1, \ldots, \sS_r\}$ 
of $K_c(X)$ which is dual to the first one. Basically, $\sS_k$ is the
class corresponding to the complex
$$
  \sR_{k}^{*} \to \bigoplus_{l} \, b_{kl} \, \sR_{l}^{*} 
              \to \bigoplus_{l} \, a_{kl} \, \sR_{l}^{*} 
  \to \sR_{k}^{*}.
$$
The Fourier-Mukai transformation induces the isomorphisms:
\begin{itemize}
 \item $\Phi: K (X) \to K^G (\IC^3)$, where $K^G (\IC^3)$ is the
 Grothendieck group of $G$-equivariant sheaves on $\IC^3$.
 It acts by $\Phi (\sR_k) = R_k \tensor \sO_{\IC^3}$.
 \item $\Phi_c: K_{c} (X) \to K^G_c (\IC^3)$, where $K^G (\IC^3)$ is
 the Grothendieck group of $G$-equivariant sheaves supported at the origin.
 It acts by $\Phi_c (\sS_k) = R_k \tensor \sO_{0}$.
\end{itemize}
Moreover each of the groups $K^G (\IC^3)$  and $K^G_c (\IC^3)$ is
isomorphic to $R(G)$, the representation ring of $G$. The following
diagram illustrates all these isomorphisms:
\medskip
\begin{equation*}
 \xymatrix{
   \sS_{R} \ar@{|->}[d] & K_c (X) \ar[d]_{\Phi_c} \ar[rr] && K(X)
   \ar[d]^{\Phi} & \sR \ar@{|->}[d] \\
   R\tensor \sO_{0} \ar@{|->}[d] &  K_c^{G} (\IC^3) \ar[d]_{\iso} \ar[rr]&&
       K^G (\IC^3) \ar[d]^{\iso} & R \tensor \sO_{\IC^3}
   \ar@{|->}[d]\\ 
   R & \  R (G) \ar[rr]^{\kappa} && R(G) & R 
         } 
\end{equation*}
\medskip
The bottom morphism $\kappa$ is the multiplication by the element 
$$
  \lambda = \sum_{i=0}^3 \, (-1)^i \, \Lambda^i Q
$$
of $R(G)$. Since $G \subset SL(3, \IC)$ we have $\lambda = - Q + \Lambda^2 Q$.
If we consider the basis $\{R_{\mbox{reg}}, R_1, \ldots, R_r\}$ on the
left-hand-side and the basis $\{R_0, R_1, \ldots, R_r\}$ on the
right-hand-side, then 
\begin{equation}\label{eq:kappa}
 \kappa = \left[
            \begin{array}{ll}
	      0 & 0\\
	      0 & C
            \end{array}
          \right].
\end{equation}

At the level of $K$-theory the bases $\{\sS_k\}$ and $\{\sR_k\}$ are
dual to each other. Therefore the corresponding $\{\ch_c(\sS_k)\} \subset H^{*}_c (X, \IQ)$
and $\{ \ch (\sR_k) \} \subset H^{*} (X, \IQ)$ are dual to each other
with respect to the pairing
\begin{equation}\label{eq:pairing}
  \< \ch_c(\sS_k), \ch (\sR_l) \> = \int_X \ch_c(\sS_k) \cup \ch(\sR_l) \, \,  \hat{A} (X).
\end{equation}
At the level at the representation ring it should be that the above
pairing descends to the multiplication of the corresponding virtual 
characters. However, the bases we gave for $R(G)$ are not dual to each other since
$\< R_{\mbox{reg}}, R_k \> = \dim R_k$ for all $k$'s. Inspired by the
multiplication formula~\eqref{eq:mult-formula} we modify the basis on
the right-hand-side of $\kappa$ to 
$$
 \{
    R_0, R_1 - \dim(R_1)\, R_0, 
   \ldots, R_r -\dim(R_r) \, R_0
   \, 
 \}.
$$
This new basis does not modify the form~\eqref{eq:kappa} for the
matrix of $\kappa$. Also, its lifting in $K(X)$ is
$$
  \{\ch(\sR_0), \ch(\sR_1) - \rk(\sR_1)), \ldots, \ch(\sR_r) -
  \rk(\sR_r)\}.
$$
Because of the properties of the Cartan matrix, this basis is still
dual to the basis of $K_c(X)$ given by $\sS_{k}$'s. Therefore we
modified the basis of $K(X)$ so that the multiplicative pairing in
$K$-theory is compatible to the natural pairing in the representation
ring.

\vskip 20pt
%\newpage
\bibliographystyle{alpha}
%\begin{thebibliography}{99}

%%%%
\end{document}